\documentclass[conference]{IEEEtran}
\usepackage[colorlinks=true]{hyperref}
\usepackage{url}
\usepackage{subfig}
\usepackage[dvips]{graphicx}
\usepackage[english]{babel}
\usepackage[latin1]{inputenc}

\usepackage{amssymb}
\usepackage{url}
\usepackage{yhmath}

\usepackage{amsmath}
\usepackage{amsfonts}
\usepackage{algorithmic}
\usepackage[ruled]{algorithm2e}
\usepackage{times}
\usepackage{color}
\usepackage{afterpage}
\usepackage{placeins}

\usepackage{cite} 
\usepackage{stfloats}  
                        
\usepackage{color}
\usepackage{amsfonts}

\usepackage{caption}
\newcommand {\bX}{\mathbb{X}}
\newcommand {\bA}{\mathbb{A}}
\newcommand {\bK}{\mathbb{K}}
\begin{document}


\title{Quantization and Stochastic Control of Trajectories of Underwater Vehicle in 
Bearings-only Tracking} 
\author{
\IEEEauthorblockN{Huilong Zhang}
\IEEEauthorblockA{Univ. Bordeaux\\
CNRS IMB, UMR 5251 \\
INRIA Bordeaux-Sud Ouest, France\\
huilong.zhang@math.u-bordeaux.fr}\\
\IEEEauthorblockN{Dann Laneuville}
\IEEEauthorblockA{DCNS Research \\
Nantes, France }
\and 
\IEEEauthorblockN{Beno\^{\i}te de Saporta}
\IEEEauthorblockA{Univ. Montpellier \\
CNRS IMAG, UMR 5149\\
INRIA Bordeaux-Sud Ouest, France\\ \\
}\\
\IEEEauthorblockN{Adrien N\`egre}
\IEEEauthorblockA{DCNS Research \\
Toulon, France }
\and
\IEEEauthorblockN{Fran\c cois Dufour}
\IEEEauthorblockA{Inst. Polytech. Bordeaux\\
CNRS IMB, UMR 5251 \\
INRIA Bordeaux-Sud Ouest, France}\\
}

%

\maketitle

\selectlanguage{english}

\begin{abstract}
We present in this paper a numerical method which computes the optimal trajectory of a 
underwater vehicle subject to some mission objectives. The method is applied to a submarine whose 
goal is to best detect one or several targets,  or/and to minimize its own detection range 
perceived by the other targets. The signal considered is acoustic propagation attenuation.  Our approach is based on 
dynamic programming of a finite horizon Markov 
decision process. A quantization method is applied to fully discretize 
the problem and allows a numerically tractable solution. Different scenarios 
are considered. We suppose at first that the position and the velocity of the targets are 
known and in the second we suppose that they are unknown and estimated by a 
Kalman type filter in a context of bearings-only tracking.
\end{abstract} 


\IEEEpeerreviewmaketitle

\section{Introduction}

Target tracking, in particular submarine target tracking, which role is to determine the 
position and velocity of the target has been extensively studied in the 
past decade \cite{ristic2004a,laneuville08a, baser2010a}. To our knowledge, there is little work \cite{negre12a} that focuses on computation of
optimal trajectories of underwater vehicle based on signal 
attenuation due to acoustic propagation and taking into account the uncertainties on the target position and velocity. The original aspect of our approach 
is to propose a mathematical computational model to address this problem. 

In a context of passive underwater acoustic warfare, we are interested in optimizing 
the trajectory of a submarine called \emph{carrier}. 
The carrier is equipped with passive sonars. 
 A question naturally arises: how should the carrier position itself 
to detect at best the acoustic signal issued by other vehicles called \emph{targets}? 
Conversely, given that the targets are
also equipped with sonars, how should the carrier position itself to keep its 
own detection range as low as possible with respect to those targets? 
These two objectives being seemingly 
contradictory, is it possible to take into account both of them simultaneously?  
A smart operator, if provided information about a single target position and velocity 
and a sound propagation code can find a good trajectory for either one of these single objectives. 
If the two criterions are considered simultaneously, or if several targets have to be 
taken into account, it is hardly possible for a human operator to find the best 
route. Some decision making tools have to be developed to this aim. 

One of the candidates are Markov decision processes (MDPs) which are widely used in many fields such as 
engineering, computer science, operation research. 
They constitute a general family of controlled stochastic 
processes suitable for modeling sequential decision-making problems under uncertainty.  
A significant list of references on discrete-time MDPs may be found in the survey \cite{arapostathis93a} and the books 
\cite{bertsekas78a, hernandez96a, hernandez99a, puterman05a,chang2013a}.
 The objective of this paper is to use this framework to compute 
optimal trajectories for underwater vehicles evolving in a given 
environment to meet some objectives such as mentioned above. 
In our preliminary works \cite{negre12a} and \cite{zhang14b} we have applied this method to 
control a carrier whose mission was to 
detect one or several targets as well as possible. 
In \cite{negre12a} one important assumption was that the positions of the targets followed a 
random model but were perfectly observed at the decision dates.
This assumption is unrealistic as the targets positions are usually only known up to some 
random error through sonar measurements. The assumption is abandoned in \cite{zhang14b}. In the present 
paper we propose a control strategy for a carrier submarine in a context of bearings-only 
tracking (BOT). 

It is well known that many real-world problems modeled by MDPs have finite but huge state and/or action spaces, 
leading to the well-known curse of dimensionality, which makes the solution of the resulting models numerically
intractable. Our initial optimization procedure includes an optimized dynamic discretization of 
the space of 
the targets positions at each time step from the current time to the computation horizon. 
Adding uncertainty on the targets makes the space of the targets positions significantly larger. 
More points are thus required in the discretization grids to keep a satisfying accuracy, dramatically 
slowing down the global process: construction and optimization of the discretized grids, final 
trajectory optimization procedure. Simulation-based algorithms such as multi-stage adaptive sampling, 
evolutionary policy method and model reference adaptive search \cite{chang2013a} can be 
used to overcome this difficulty. But in the context of BOT, these approaches can not be applied, 
since a filtering procedure has to be integrated. Under the consideration of practice, 
we choose to split the long-term horizon into a sequence of sub-intervals, making the 
global optimization sub-optimal but real-time achievable, as the discretization grids only need to 
be optimized on a shorter horizon.

Related works in the literature are for example 
model predictive control \cite{camacho2013a}, moving horizon control \cite{vandenbroeks2002a} or 
receding horizon approach \cite{chang2003a}, where similar ideas of shortened horizon are used to make the problems numerically solvable in real-time. within the POMDP (partially observable Markov decision processes 
\cite{bertsekas2005a, lovejoy1991a}) framework, the most computationally tractable strategy is to use myopic open-loop 
feedback control. In \cite{beard2015a}, a sequence of shortened horizon H is fixed and Monte Carlo integration method is used 
to estimate reward fonction for each horizon. 

Our focus is on computational approaches to find optimal policies. The paper is organized as follows. 
Section \ref{sect_statement} states the problem in a general context and explains our strategy to 
solve it. Section \ref{sect_controle_stoch} briefly outlines the dynamic programming algorithms used. 
Section \ref{sect_numerical_results} presents the numerical results for four 
different scenarios, from the simplest case to the realistic case.
\section{Problem Statement}
\label{sect_statement}
\begin{figure} [t]
\begin{center} 
\includegraphics[width=9cm]{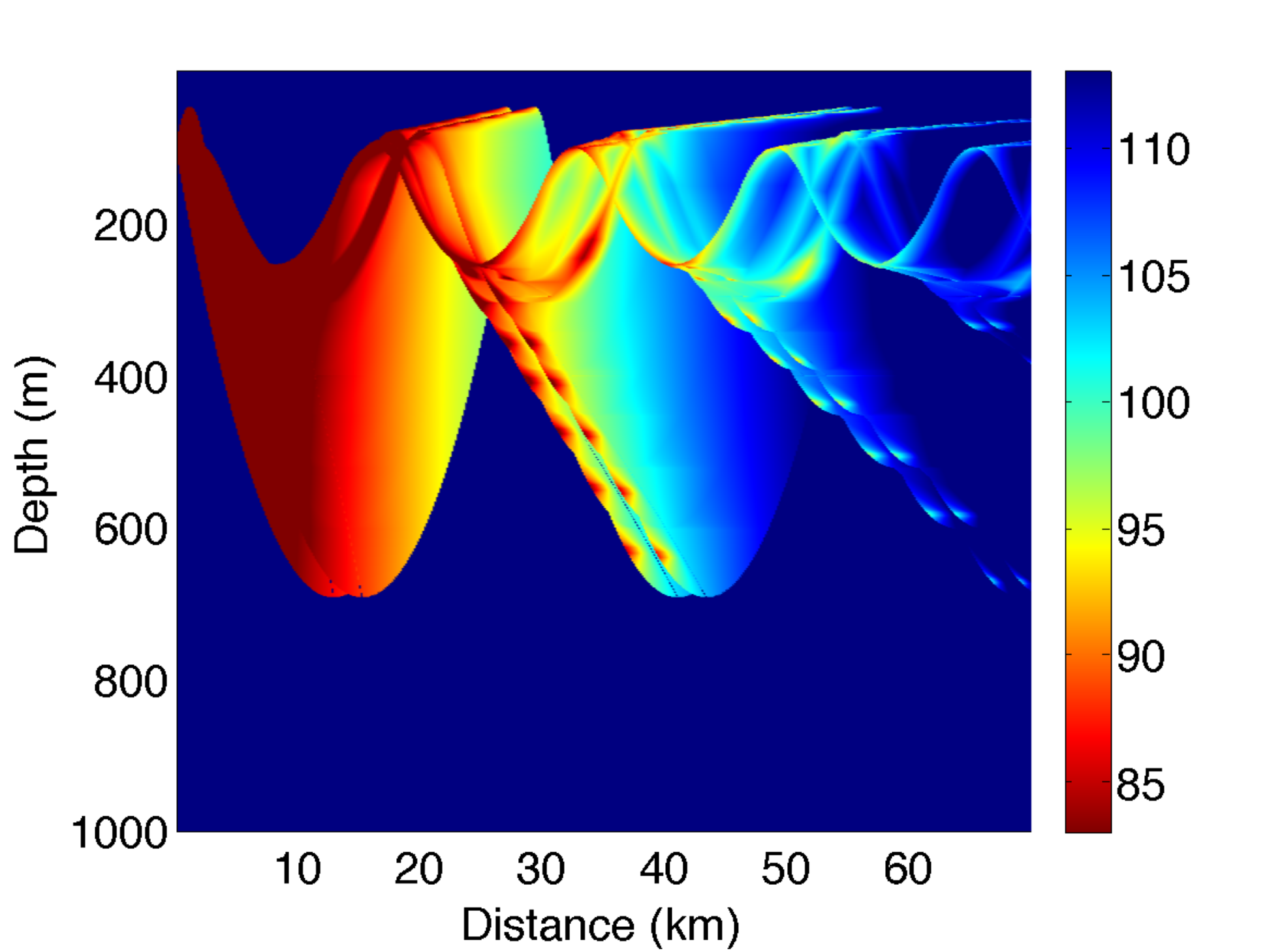}
\caption{Sound propagation diagram: signal loss level vs depth \& distance}
\label{fig_1}
\end{center} 
\end{figure}

Underwater sound propagation in the sea depends on many environmental parameters such as 
temperature, salinity, etc. Indeed, small variations of these parameters may greatly modify the 
sound propagation (see for instance \cite{etter2013a}). To illustrate this phenomena, 
Fig.~\ref{fig_1} gives an example of a sound propagation diagram with an emitter 
source standing at a $100\,m$ depth (top left). 
The signal loss level due to propagation usually takes values in $[80,200]\,dB$.
The emitter may well be detected by a carrier in the areas (distance to emitter versus depth of carrier) of lower signal loss level 
($\approx 80\,dB$, represented by dark red color), and a contrario is hardly detectable if not undetectable 
in higher signal loss level areas ($\approx 200\,dB$, represented by dark blue color). Please also note that the carrier being close to the emitter does not guarantee  a good detection. In this figure, the signal 
loss level is saturated to $120\,dB$.

Let us consider a general situation with a carrier submarine of interest $S$ 
surrounded by one or several targets. This submarine carries passive sensors so that the 
available measurements are the estimated angles (bearings) from observer to sources. 
The trajectory of submarine $S$ has to be controlled in order to satisfy some given mission 
objectives. These can be e.g. optimizing the different targets
detection range, minimizing its own detection range as perceived by the other 
targets, reaching a way-point with minimum fuel consumption, etc. 

The following inputs are supposed to be available
\begin{enumerate} 
  \item $S$ perfect position and velocity,
  \item noisy observation of the target, namely bearing and frequency,
  \item information about the environment (sound speed, sea floor depth, bottom type, ...).
\end{enumerate} 
A sound propagation code which calculates the signal loss level due to propagation 
(such as in Fig.~\ref{fig_1}) is also available. 

The desired output is an optimized sequence of commands to pilot the carrier submarine in 
order fulfill its mission at best within a given time interval, see Fig 
\ref{fig_2}.  

Concerning the targets, their evolution is described by a stochastic model. 

\section{Stochastic optimal control}
\label{sect_controle_stoch}
As seen in Section \ref{sect_statement} the problem is quite complex, there are 
uncertain inputs and a probabilistic modeling should be used to take account of this 
feature. A discrete stochastic optimal control framework seems to be a good
candidate to solve this problem. More precisely, we use a discrete-time finite horizon 
dynamic programming approach. We first describe this approach for a single long-term optimization problem, then we will explain how to divide it into a sequence of short-term ones.

\begin{figure} [t]
\begin{center} 
 \includegraphics[width=9cm]{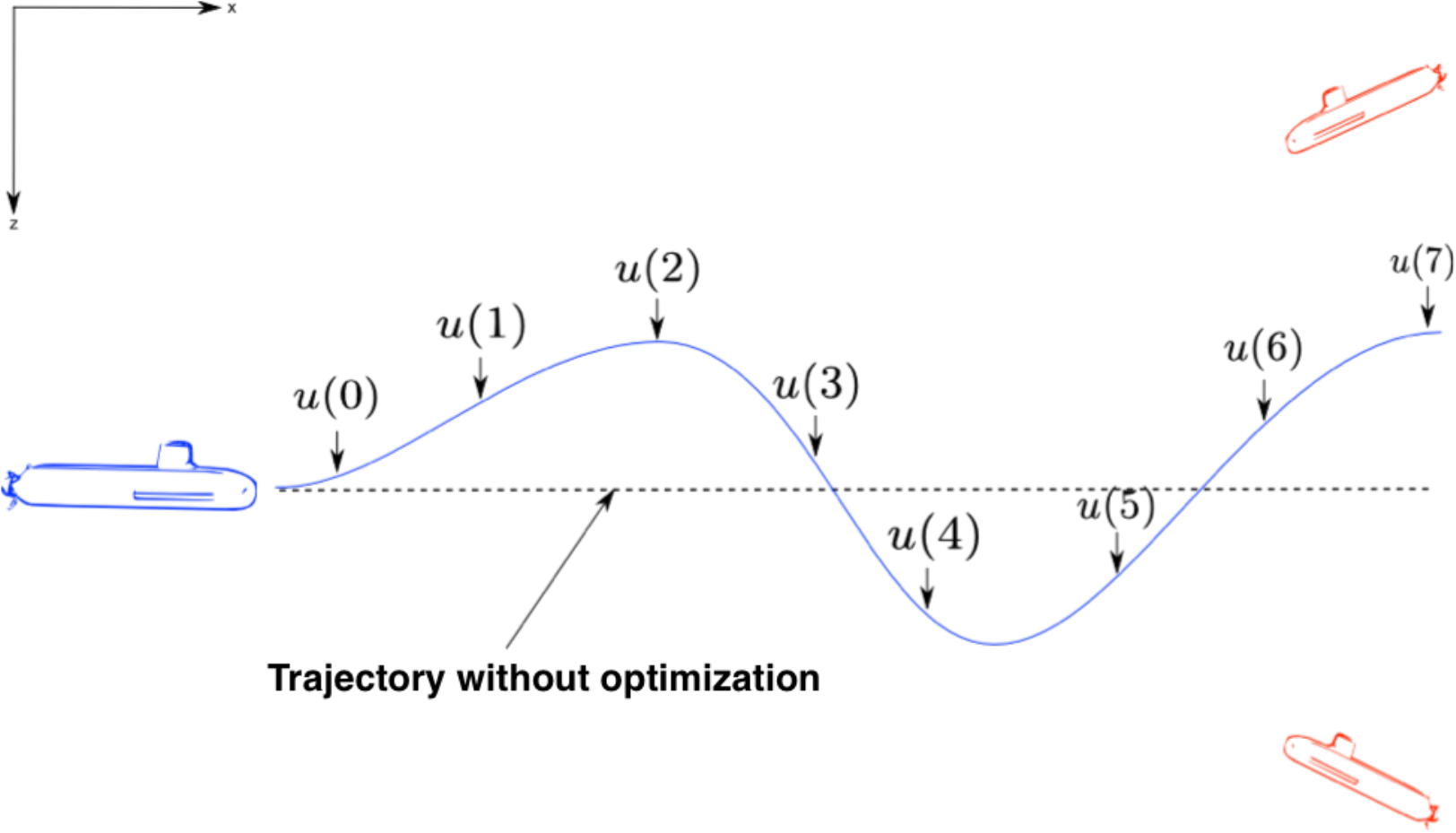}
\caption{Optimal trajectory example}
\label{fig_2}
\end{center} 
\end{figure}

We have opted for the quantization technique to 
discretize the movement of the targets. As the targets evolve in an infinite 
continuous space, quantization offers both a discretized density and a transition 
matrix, required for MDP approach and allows for discretization grids that are 
dynamic and concentrated at the most probable positions of the targets. The aim is to 
compute an optimal trajectory by applying a command $u(t)$ to $S$'s future states for 
all times $t$ up to the computation horizon, see Fig. \ref{fig_2}.

In our framework, three probabilistic techniques are used, namely Kalman filter, 
quantization and dynamic programming, and each of these overlaps with 
one another and are combined in an iterative algorithm. 

\subsection{Dynamic programming}\label{sec:DP}
In this section, we briefly introduce the discrete-time finite horizon Markov 
control framework. Let us consider the following model: 
\begin{displaymath} 
  (\bX, \bA, \{A(x)|x\in \bX\}, P, c)
\end{displaymath} 
with 
\begin{itemize} 
\item $\bX$: finite space, namely the state space,
\item $\bA$: finite space representing the control or actions set,
\item $\{A(x)|x\in \bX\}$: family of non empty subsets of $\bA$, $A(x)$ being the set of feasible controls or actions when the system is in state $x\in \bX$. We suppose that $\bK : = \{(x,a)|x\in \bX, a \in A(x)\} $ is a measurable subset of  $\bX\times \bA$,
\item $P$: transition matrix on $\bX$ given $\bK$ which stands for the transition probability function,
\begin{equation}\label{eq_kernel}
P_{xy}^t(a) = \mathbf{P}(x_{t+1} = y |x_t = x, a)
\end{equation}
\item $c$: $\bK \rightarrow \mathbf{R}$ measurable function representing the cost per stage,
\item $C_N$: $\bX \rightarrow \mathbf{R}$ measurable function representing the terminal cost.
\end{itemize} 

Suppose that a finite horizon $N\geq 1$ and an initial state $x\in \bX$ are given. 
The total expected cost of a policy $\pi = \{a_t|t=0,1,\cdots,N-1\} \in \Pi$ 
($\Pi$ is the set of all possible policies) is defined as 
\begin{equation} 
\label{cout_moyen}
\begin{array} {rcl}
\displaystyle 
J(\pi,x) & = & \mathbf{E}_{x}^{\pi} \left [\displaystyle \sum_{t=0}^{N-1} c(x_t,a_t)+C_N(x_N)\right ] \\
 & = & \mathbf{E}\left [\displaystyle \sum_{t=0}^{N-1} c(x_t,a_t)+C_N(x_N)| \pi, x_0=x \right ].
\end{array} 
\end{equation} 
The optimal total expected cost function is then defined as 
\begin{displaymath} 
  J^\ast(x) = \underset{\pi \in \Pi}{\mbox{opt}} J(\pi,x), ~\forall x \in \bX.
\end{displaymath} 
The goal is to find an optimal policy $\pi ^\ast$ which optimizes this cost. 
The optimization criteria $\mbox{opt}$ can be either $\mbox{inf}$ (to minimize a loss) or $\mbox{sup}$ (to maximize a performance) according to the mission's type.  The main 
idea of the dynamic programming principle is to recursively solve a sequence of one-step optimization problems by
computing (backwards) a sequence of functions $J_N^*, J_{N-1}^*, \cdots, J_0^*$
by the following algorithm
\begin{itemize} 
  \item $J_N^\ast(x) = C_N(x), \forall x\in \bX$,
  \item for $t\in \{N-1,\cdots,0\}$ and $x\in \bX$
\begin{equation} 
\label{eq_retrograde}
  J_t^\ast (x) = \underset{\pi \in \Pi} {\mbox{opt}} 
  \displaystyle \sum_{y\in \bX} P_{xy}^t(a)[c(x,a,y)+J_{t+1}^\ast(y)].
\end{equation} 
\end{itemize} 
The last function thus constructed, $J_0^\ast(x)$, is the value function $J^*$. 
In addition, the policy $\pi^*$ defined by
$\pi ^\ast = \{a_0^\ast, a_1^\ast,\cdots, a_N^\ast\}$ with 
\begin{equation} 
\label{eq_arg_opti}
a_t^\ast(x) = \underset{a\in A(x)} {\mbox{argopt}}
\displaystyle \sum_{y\in \bX} P_{xy}^t(a)[c(x,a,y)+J_{t+1}^\ast(y)]
\end{equation} 
is an optimal policy, i.e. $J(\pi^*,x)=J^*(x)$.
\subsection{Submarine $S$}\label{sec:subm}
The movement of the carrier is assumed to be deterministic, modeled by the $3$-dimensional vector $s_t$ for the carrier position at time $t$. Its dynamic is determined by the control policy
\begin{displaymath} 
  s_{t+1} = s_t + a_t, ~s_0 ~\mbox{known}. 
\end{displaymath} 
At each time step and at each position of sensor $s_t$ the control set is
\begin{eqnarray*} 
\lefteqn{\bA=\{(i_1\Delta_{s_1},i_2\Delta_{s_2},i_3\Delta_{s_3}),}\\
&&  i_j \in [-L_j, L_j-1,\ldots,-1,0,1,\ldots,L_j], j=1,2,3\},
\end{eqnarray*} 
for some given space steps $\Delta_{s_j}$ (related to the maneuverability of the carrier in the prescribed direction during the computation time step) and control ranges $L_j$.
The feasible controls $a$ at position $s$ are those that correspond to \emph{actual} underwater positions $s+a$ (e.g. not above the sea level, not below the sea floor, not on shore, \ldots).
Therefore, we can consider that the submarine $S$'s state takes values at time $t$ in a $3$-dimensional
grid ${X}^S_t$ whose size increases gradually 
\begin{displaymath} 
  {X}^S_t = \{s_0 + (n_1 \Delta_{s_1}, n_2 \Delta_{s_2}, n_3 \Delta_{s_3}), 
n_j \in [-tL_j, tL_j]\}.
\end{displaymath} 

\subsection{Targets}
The modeling of the targets being similar and independent of each other, we will consider only one target in this section. 
The multi-target case will be discussed in Section \ref{subsect_scenario_2}.
We suppose that the target evolves with a constant depth and that its 
dynamics is independent of the carrier $S$. 
The target kinematics are modelled using 2-dimensional Cartesian position and velocity vectors 
$w_t = (w_{x,t}, \dot{w}_{x,t}, w_{y,t},,\dot{w}_{y,t})^\prime$. Its evolution is described by the stochastic model
\begin{equation} 
\label{eq_target}
  w_{t+1} = F w_t + K \varepsilon _t,
\end{equation} 
where $F = 
\left [\begin{array} {cc} 
1 & T \\
0 & 1
\end{array} 
\right ]
\otimes I_2, ~~
K = \left [
\begin{array} {c} 
T^2/2 \\ T 
\end{array}
\right ] \otimes I_2$, 
$w_0 \sim {\mathcal N}(\mu_0,\sigma_0)$, $\varepsilon _t \sim {\mathcal N}(0, \Sigma_\varepsilon )$ is a 2-dimensional independant 
noise with $\Sigma_\varepsilon = \sigma_\varepsilon ^2I_2$. The parameters 
$\mu_0, \Sigma_0, \sigma_\varepsilon $ are given a priori 
standard  values.

Thus the natural state space for the target is continuous.
However, to apply the dynamic programming approach described in part~\ref{sec:DP}, we need a finite state space. Hence, it is necessary to 
approximate the continuous state space with a finite one.

Various methods exist for this approximation \cite{kushner92a}. We chose a quantization 
method as it is dynamic and suitable to discretize random processes.  
The goal of this method is to approximate the continuous 
state space Markov chain $w_t = (w_{x,t},\dot{w}_{x,t},w_{y,t},\dot{w}_{y,t})^\prime $ by a finite state space chain
$\hat{w}_t = (\widehat{w}_{x,t}, \widehat{\dot{w}}_{x,t}, \widehat{w}_{y,t},  \widehat{\dot{w}}_{y,t})$. To this aim we use the quantization 
algorithm described in \cite{pages04a,pages04b, pages05a}. Roughly speaking, 
more points are put in the areas of high density of the random variable. The 
quantization algorithm is based on Monte Carlo simulations combined with a 
stochastic gradient method. It provides $N+1$ grids, one for each $w_t$ 
($0\leq t \leq N$), with a fixed finite number of points in each grid. The algorithm 
also provides weights for the grid points and probability transition between 
two points of two consecutive grids, thus fully determining the distribution of 
the approximating sequence $\hat{w}_t$ ($0\leq t \leq N$). The quantization theory ensures that 
the $L^2$ distance between $\hat{w}_t$ and $w_t$ tends to 0 as the number of 
points in the quantization grid tends to infinity \cite{pages04b}. In our case, 
the process $w_t$ is approximated by a finite Markov chain with the following notation
\begin{itemize} 
  \item $\hat{w}_t= (\widehat{w}_{x,t}, \widehat{\dot{w}}_{x,t}, \widehat{w}_{y,t}, \widehat{\dot{w}}_{y,t})^\prime$ is the quantization of $w_t$ and $\hat{X}^W_t$ is the 
$M$ point grid at time $t$,
  \item $\hat{P}^t$ is the transition matrix, 
$\hat{P}^t_{ij} = \mathbf{P}(\hat{w}_{t+1} = j | \hat{w}_{t} = i)$, 
$\forall (i,j) \in \hat{X}^W_t\times \hat{X}^W_{t+1}$.
\end{itemize} 
There exists extensive literature on quantization methods for random variables and processes \cite{gray1998a,graf2000a,pages04a,pages04b, pages05a,saporta2015a}. 

\subsection{Markov control model}
We can now fully describe the Markov control model corresponding to our problem.
The finite state space is $\bX =\cup_{t=0}^N\bX_t$, where $\bX_t={X}_t^S\times 
\hat{X}^W_t$.
The process at time $t$ is $x_t := (s_t, \hat{w}_t)$. Component $s_t$ evolves deterministically while $\hat{w}_t$ is 
stochastic.  
The action space is $\bA$ described in Section~\ref{sec:subm}.

The process $\hat{w}_t$ is not controlled hence the transition matrix $P$ can 
be obtained from $\hat{P}$ as
\begin{eqnarray*}
P_{xy}^t(a)&=&\mathbf{P} (x_{t+1}=y|x_{t}=x, a_t=a)\\
&=&1_{\{m=l+a\}}\hat{P}^t_{ij},
\end{eqnarray*}
for all $x=(l,i)\in\bX_t$ and $y=(m,j)\in\bX_{t+1}$.

The transition cost function from state $x = (l,i)$ to state 
$y=(m,j)$ given a control $a$ is $c(x,a,y)$.
In our application, this function only depends 
on the sound propagation code  which is \emph{geometric} (its values only depend 
on the submarine and target relative positions and some fixed parameters). Hence $c(x,a,y) = c(m,j)$. We
finally obtain the following dynamic programming equation
\begin{equation} 
\label{equ_belleman}
J_t^\ast(l,i) =  \underset{a\in \bA(l,i)} {\mbox{argopt}}
  \displaystyle \underset{j\in \hat{X}^W_{t+1}} {\displaystyle \sum} 
\hat{P}^t_{ij}[c(l+a,j)+J_{t+1}^\ast(l+a,j)].
\end{equation} 

\subsection{Algorithm}
\label{subsect_algorithm}
In the case where the position and velocity of the target are known, the 
resolution of our stochastic optimal control problem can be divided in three steps. 
The first step is an optimal quantization in order to approximate the target state by a finite 
Markov chain state. The second step is a backward dynamic programming algorithm for evaluating the 
best policy from each possible system state. This step only depends on the law of the process.
The last step is on line, it consists in applying the optimal control policy to a given trajectory of 
the process. 

The case where the state of the target is unknown is described in \ref{subsect_scenario_3} and \ref{subsect_scenario_4}.

\subsection*{Step 1: Quantization}
The goal is to approximate $w_t$ by a finite Markov chain $\hat{w}_t$. Algorithm~\ref{algo_quantization} 
describes how to adapt the CLVQ (competitive learning vector quantization) technique to the context of an 
optimal quantization of a Markov chain. We refer the reader to \cite[section 3.3]{bally03a} for the computation 
of the weights and transition matrices by Monte Carlo simulations.
\begin{algorithm}[ht]
\SetKwInOut{Input}{input}
\SetKwInOut{Output}{output}
\SetKwInOut{Return}{return}
\caption{Extended CLVQ algorithm\index{Extended CLVQ algorithm} to quantize a Markov chain}
\label{algo_quantization}
\Input{Number of point $M$, Number of runs $NR$, Horizon $N$, Sequence $(\gamma_t)$\\
Initial grids $(\hat{X}^W_t)$ $0\leq t\leq N$ with $M$ points, \\
Simulator of trajectories of target Markov chain}
\Output{Optimized grids $(\hat{X}^W_t)$, $0\leq t\leq N$}
\Begin{
\For{$m \leftarrow 0$ \KwTo $NR-1$}{
simulate trajectory $(w_0,w_1,\ldots,w_N)$ according to law of Markov chain\\
\For{$t\leftarrow 0$ \KwTo $N$}{
\emph{competitive phase} select $y$ closest neighbour of $w_t$ in $\hat{X}^W_t$\\
\emph{learning phase} set $y'=y-\gamma_{t}(y-w_t)$\\
$\hat{X}^W_{t} \leftarrow \hat{X}^W_t\cup\{y'\}\backslash\{y\}$
}
}
\Return{$(\hat{X}^W_t)$, $0\leq t\leq N$}
}
\end{algorithm}

We obtain for each time 
$t=0,1,\cdots, N$ a grid $\hat{X}^W_t$ and the corresponding transition matrix $\hat{P}^{t}_{ij}$. Grid 
and transition matrices are stored. Let $\hat{X}^W_{t}$ be the optimal $L_{p}$-quantization of the random variable $w_{t}$ 
for $0\leq t\leq N$ by a random variable taking $M$ points.
Let us denote by $\hat{X}^W_{t}=\{\hat{w}_{t}^{1},\ldots,\hat{w}_{t}^{M}\}$ the grid at step $t$ and the associated Voronoi 
tessellation $C_{1}(\hat{X}^W_{t}),\ldots,C_{M}(\hat{X}^W_{t})$. Clearly, the process $(\hat{X}^W_{t})_{0\leq t\leq N}$ is 
not a Markov chain. However, at each step $0\leq t\leq N-1$ one can compute
\begin{align*}
\hat{P}^{t}_{i,j} 
& =  \mathbf{P}\big(\hat{X}^W_{t+1}=j \big| \hat{X}^W_{t}=i \big) 
=\frac{\mathbf{P}\big(\hat{X}_{t}=i, \hat{X}_{t+1}=j \big)}{\mathbf{P}\big(\hat{X}_{t}=i\big)}  \\
& =  \frac{\mathbf{P}\big(w_{t}\in C_{i}(\hat{X}^W_{t}), w_{t+1}\in C_{j}(\hat{X}^W_{t+1}) \big)}{\mathbf{P}\big(w_{t}\in C_{i}(\hat{X}^W_{t})\big)},
\end{align*}
for $1\leq i,j\leq M$.
The marginal quantization approximation of the Markov chain $(w_{t})_{0\leq n\leq N}$ is then defined by the Markov 
chain $(\hat{X}^W_{t})_{0\leq n\leq N}$ whose transition matrix at step $t$ is given by
$\big (\hat{P}^{t}_{i,j}\big)_{0\leq t \leq N-1}$, 
and the initial distribution of ${w}_{0}$ is given by that of $\hat{X}^W_{0}$. It can be shown (see for all the details \cite[Theorem 3.1]{pages04b}) that the joint distribution of $(w_{t})_{0\leq t\leq N}$ can be approximated by that of 
$(\hat{X}^W_{t})_{0\leq t\leq N}$ with a rate of convergence of order $N^{1+{1}/{d}}M^{-1/d}$, where $d$ is the dimension of the state space.

\subsection*{Step 2: Pre-computations}
Compute for each possible state element 
$$ (l,i) \in \bX_t = X_t^S\times \hat{X}_t^W, ~0\leq t\leq N$$ the optimal control policy that optimizes  the cost function
for horizon $N$.
\begin{algorithm}[ht]
\SetKwInOut{Input}{input}
\SetKwInOut{Output}{output}
\SetKwInOut{Return}{return}
\caption{Dynamic programming} 
\label{algo_programmation_dynamique}
\Input{Number of point M, Horizon $N$, \\
Grids $(\hat{X}^W_t)_{0\leq t\leq N}$, \\
Transition matrix $(\hat{P}^t_{ij})_{0\leq t\leq N}$}
\Output{Optimal control $a_t^\ast$}
\Begin{
\For{$(l,i)\in \bX_N$} {
calculate final cost function $C_N(l,i)$ 
\begin{displaymath} 
\begin{array} {l}
\displaystyle 
  J_N^\ast(l,i)  =  \displaystyle C_N(l,i),~~\forall (l,i)\in \bX_N\\
  a_N^\ast(l,i)  =  0  
\end{array} 
\end{displaymath} 
}
\For {$t= N-1$ \KwTo 0}{ 
Calculate for $ (l,i)\in \bX_t$
{\small
\begin{displaymath} 
\begin{array} {l} 
\displaystyle 
      J_t^\ast(l,i)  = \underset{a\in \bA(l,i)}{\hbox{opt}}
      \sum_{j \in \hat{X}^W_{t+1}}\hat{P}_{ij}^t 
\big [c(l+a,j) + J_{t+1}^\ast(l+a,j)\big ]\\
      a_{t}^{\ast}(l,i)  =  \displaystyle \underset{a \in \bA(l,i)}{\mbox{argopt}}
\sum_{j \in \hat{X}^W_{t+1}}
\hat{P}_{ij}^t \big [c(l+a,j) + J_{t+1}^\ast(l+a,j)\big ]
\end{array} 
\end{displaymath} 
}
}
\Return{$a_t^\ast(l,i), 0\leq t \leq N, ~(l,i) \in \bX_t$}
}
\end{algorithm} 

\subsection*{Step 3: On line selection of optimal action}
As we suppose that the trajectory of the target is known, the selection of the optimal action given the current state and time (see algorithm 
\ref{algo_controle_optimale_en_ligne}) can be done on line. 

\begin{algorithm}[ht]
\SetKwInOut{Input}{input}
\SetKwInOut{Output}{output}
\SetKwInOut{Return}{return}
\caption{Optimal control with known target} 
\label{algo_controle_optimale_en_ligne}
\Input{Optimal control $a_t^\ast$\\
Initial condition of the system $(s_0,w_0)$\\
Grids $\hat{X}_t^W, 0\leq t \leq N$}
\Output{Adaptative control $\hat{a}_t, 0\leq t \leq N$\\
Adaptative $S$ position $s_t, 0\leq t\leq N$}
\Begin{
\For{$t=0$ \KwTo $N-1$}{
Find $\hat{w}_t$ the nearest neighbor of $w_t$ among grid $\hat{X}^W_t$\\
$\hat{a}_t = a_t^\ast(s_t,\hat{w}_t)$ optimal control at time $t$\\
$s_{t+1} = s_t+\hat{a}_t$, update $S$ position
}
}
\end{algorithm}

\section{Numerical results}
\label{sect_numerical_results}
In order to illustrate numerical results, several scenarios with different mission objectives 
have been considered. 
\begin{enumerate} 
 \item Scenario 1 : one target, its state (position and velocity) is known, 
objective is to detect at best the acoustic signal issued by the target.
\item  Scenario 2 : two targets, their states are known,  objective is to detect at best 
the acoustic signal issued by the two target.
\item  Scenario 3 : one target, its state is unknown, objective is to detect at best the acoustic signal issued by the  target.
\item  Scenario 4 : one target, its state is unknown, objective is to detect at best 
the acoustic signal issued by the  target and at same time to keep its own detection as 
low as possible.
\end{enumerate} 
The first scenario allows to validate the algorithms presented in Section \ref{subsect_algorithm}. A special 
technique of quantization is proposed for the multi-target scenario. The third and fourth scenarios are more realistic, we 
will split the time horizon into subintervals to make the computation feasible.
\subsection{Scenario 1}
\begin{figure}[t]
\begin{center} 
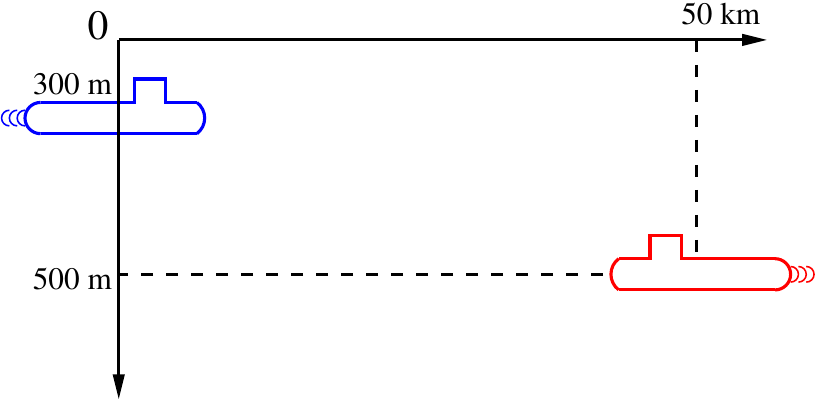
\caption{Scenario 1}
\label{fig_3}
\end{center} 
\end{figure}

We consider here the simplest scenario : there is a single target in the environment. The carrier wants to best detect 
the target. The initial geometry is depicted by Fig. \ref{fig_3}. The target (red) follows a uniform motion and its depth remains 
constant (500 m). At initial time the depth of the carrier submarine is 300 m. The relative target-carrier speed is 
10 $\mbox{m.s}^{-1}$. 

\begin{figure}[ht]
\begin{center} 
\includegraphics[width=9cm]{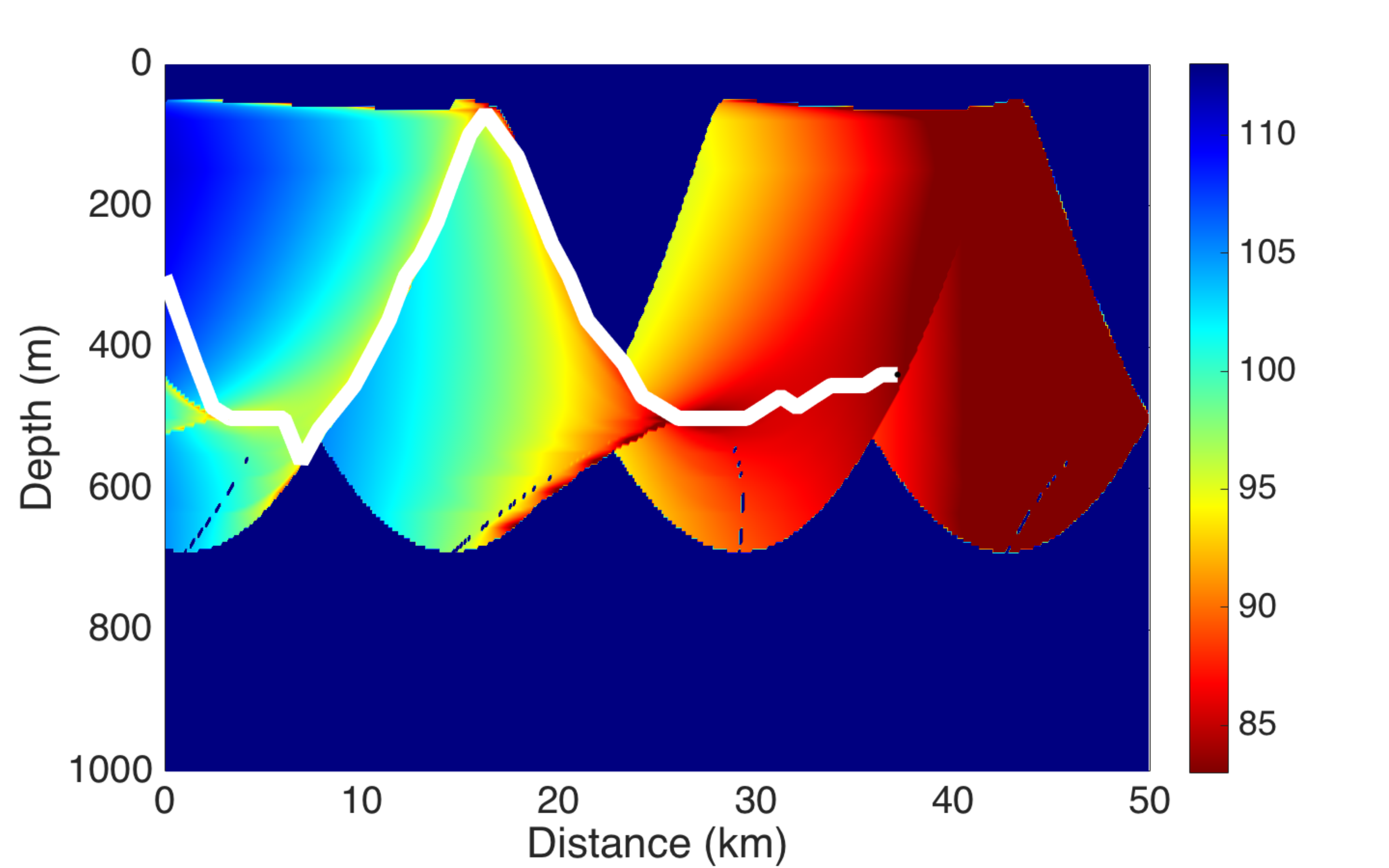}
\caption{Result of scenario 1}
\label{fig_4}
\end{center} 
\end{figure}

In this scenario we assume that the position of target $w_t$ is known at the 
decision time, and at each time step, available controls are in one 
dimension ($L_1=L_2=0$), which means that the only control action is change of depth.  
Results are illustrated in 
Fig \ref{fig_4}. The diagram represents the sound propagation 
loss of target, and the trajectory of carrier is plotted by white dot points. We can see that the movement
of carrier is ``intelligent'' in the sense that its position is always in the red aera, where the sound 
acoustic loss of target is lower.

\subsection{Scenario 2}
\label{subsect_scenario_2}
\begin{figure}[ht]
\begin{center} 
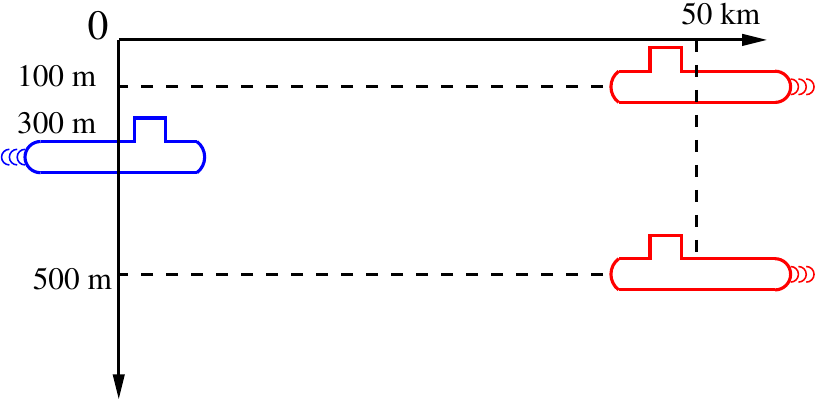
\caption{Scenario 2}
\label{fig_5}
\end{center} 
\end{figure}
Let us consider a more complex scenario with two targets at different depth. 
The first target is at 500 m depth whereas the second target has a depth of 100 m (see Fig. \ref{fig_5}).  In this scenario, the available 
controls are also in one dimension. 

In the case of multi-target, there are two ways to quantify the 
movement of targets. The first ways is to quantify each target separately and independently, because 
the evolution of the targets is assumed to be independent. The second method is to regroup the 
targets in one state space by increasing the dimension of state space. For this scenario we can 
quantify the movement of two targets ($w_t^1$ and $w_t^2$) by model
\begin{equation}
\label{eq_two_target} 
\left [
\begin{array} {c}
  w_{t+1}^1\\
  w_{t+1}^2
\end{array} 
\right ]
  = 
\left [
\begin{array} {cc}
F  &  0 \\
0  &  F 
\end{array} 
\right ]
\left [
\begin{array} {c}
  w_{t}^1\\
  w_{t}^2
\end{array} 
\right ]
+ 
\left [
\begin{array} {cc}
K  &  0 \\
0  &  K
\end{array} 
\right ]
\left [
\begin{array} {c}
  \varepsilon _{t}^1\\
  \varepsilon _{t}^2
\end{array} 
\right ]
\end{equation} 
This choice is justified by following considerations~: the first method of quantization is 
numerically intractable, since two quantization 
grids ($\widehat{X}^{W,1}_t$ and $\widehat{X}^{W,2}_t$) will be necessary and the cardinal of the
set $\widehat{X}^{W,1}_t\times\widehat{X}^{W,2}_t$ is too great to be numerically feasible. 
Using the second quantization method, some numerical dependence will be introduced and dimension of model (\ref{eq_two_target}) is twice that of the previous example. Let $\widehat{X}_{t}^{(w^1,w^2)}$ be the quantization of quantization at time $t$ calculated by model (\ref{eq_two_target}). Then a point of $\widehat{X}_{t}^{(w^1,w^2)}$ represents
a couple $(w_t^1,w_t^2)$. Let's denote by 
\begin{itemize} 
\item $i=\widehat{(w_t^1,w_t^2)} \in \widehat{X}_{t}^{(w^1,w^2)}$ a 
quantization point \`a time $t$
\item $\hat{P}^t = (\hat{P}_{ij}^t)$ the transition matrix 
\begin{displaymath} 
\displaystyle \hat{P}_{ij}^t = \mathbf{P} (\widehat{(w_{t+1}^1,w_{t+1}^2)} = j | \widehat{(w_{t}^1,w_{t}^2)} = i),
\end{displaymath} 
$\forall (i,j) \in \widehat{X}_{t}^{(w^1,w^2)}\times  \widehat{X}_{t+1}^{(w^1,w^2)}.$
\end{itemize} 
So the dynamic programming equation (\ref{equ_belleman}) remains unchanged.

The multi-target cost function is defined in order to keep a 
good detection range of each target. Let $l$ be the position of carrier, 
$i=(i_1,i_2)$ the position of targets 1 and 2, $\beta(l,i_1)$ and $\beta(l,i_2)$ 
acoustic loss value of targets 1 and 2 regarding the carrier, then the cost 
function $c(l,i)$ can be calculated by 
\begin{equation} 
\label{equat_fonct_cout_two_targets}
c(l,i) = \alpha_1 \beta(l,i_1)+\alpha_2 \beta(l,i_2), ~~~\alpha_1+\alpha_2 = 1
\end{equation} 
We take $\alpha_1=\alpha_2=\frac12$. The coefficients $\alpha_1$ and $\alpha_2$ 
can be modified if we privilege the detection one of two targets.
\begin{figure}[t]
\begin{center} 
\includegraphics[width=9cm]{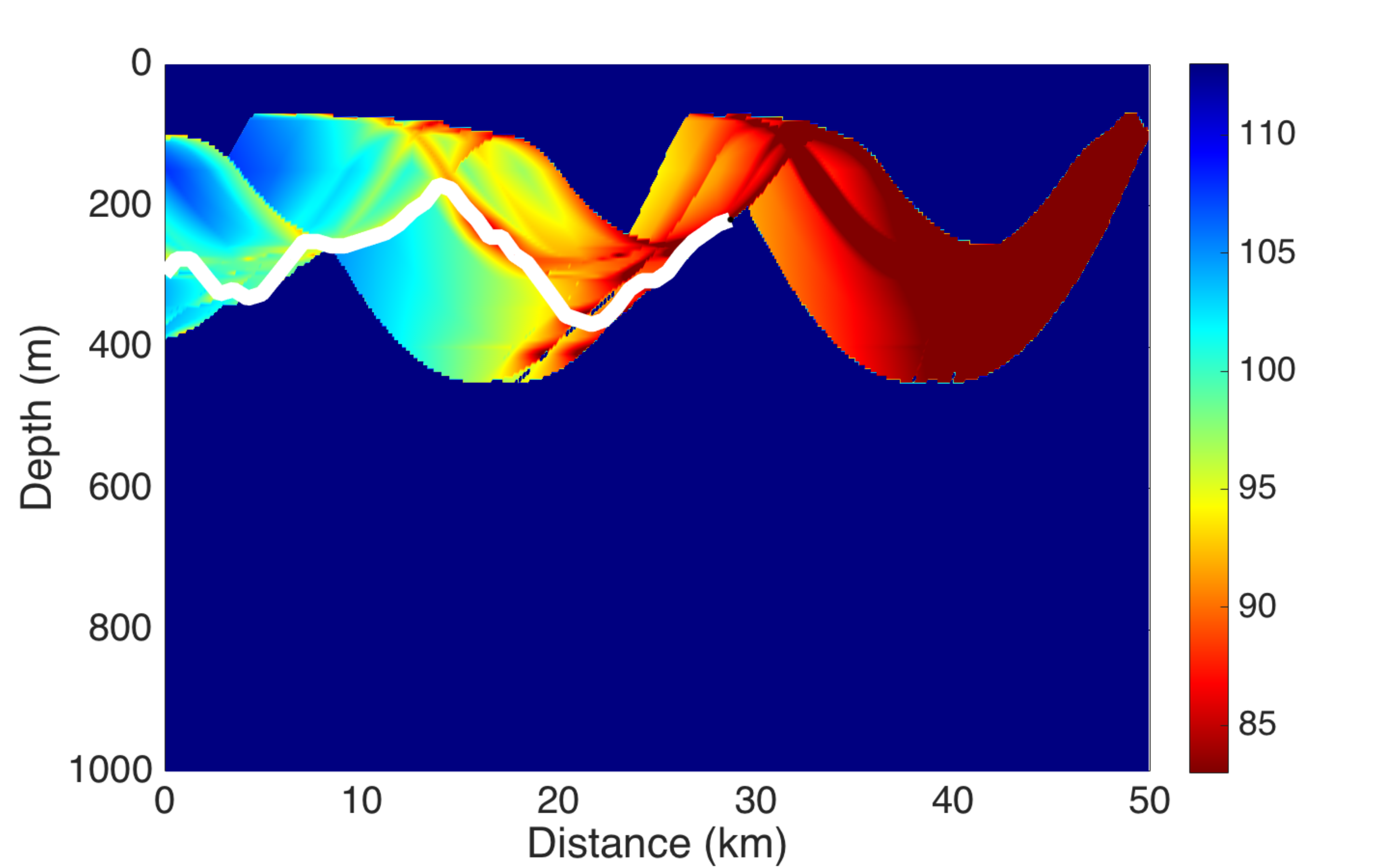}
\caption{Result of scenario 2, loss diagram of target 1}
\label{fig_6}
\end{center} 
\end{figure}

\begin{figure}[ht]
\begin{center} 
\includegraphics[width=9cm]{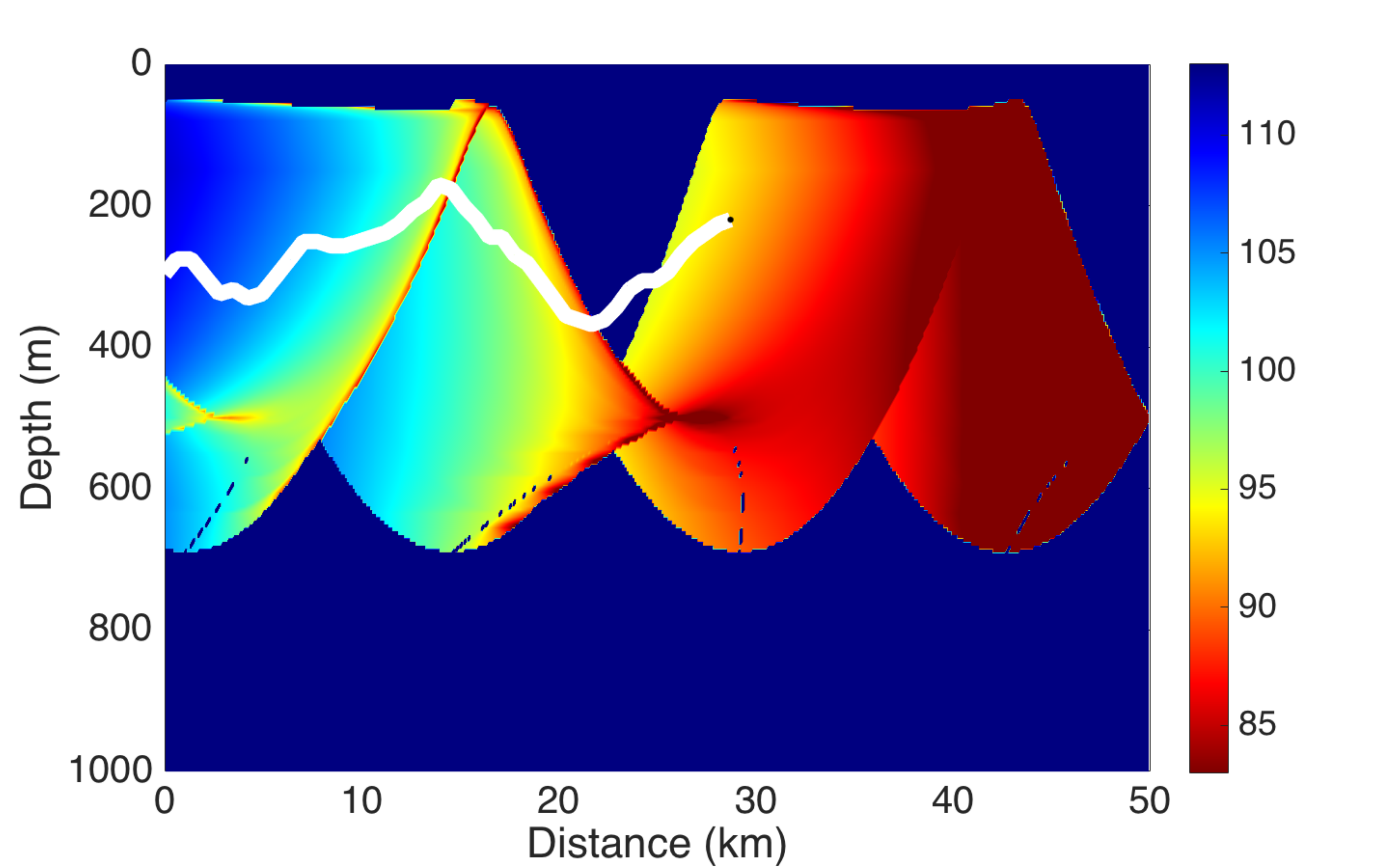}
\caption{Result of scenario 2, loss diagram of target 2}
\label{fig_7}
\end{center} 
\end{figure}
Fig. \ref{fig_6} and Fig. \ref{fig_7} illustrate the sound propagation loss diagram of targets 1 et 2 respectively, 
and the trajectory of carrier is plotted by white dot points. We can see that the movement of carrier is also ``intelligent'' in the sense that for most of the time, it is positioned in ``good'' area, where the acoustic loss is lower
 for both targets. But during certain periods, it is physically impossible to stay 
 simultaneously in red areas for both targets.  For instance in 
Fig. \ref{fig_7} when the distance with target is at 20 km, 
it is impossible to maintain the carrier in the reddish area of the two targets. However note that the carrier never reaches the blue area of both targets simultaneously.

\subsection{Scenario 3}
\label{subsect_scenario_3}
In the last two  scenarios we consider that the position and velocity of target are unknown. This assumption is more realistic because the targets positions are usually only known up to some random error through sonar measurements. 
In this article we propose a control strategy for a carrier submarine in a context of bearings-only tracking (BOT).

It is well known that when the target and carrier are under uniform rectilinear motion, the filter of BOT does not 
converge because the problem is not obervable,  some maneuvers have to be performed by the carrier in order to 
 get a precise enough approximation. Therefore we divide our operations into two periods.
In the first period there is no control, only the target motion analysis (TMA) filter algorithm is run. 
In the second period both filter and control algorithms are run, the latter being fed by the former. 

In scenario 1 and 2,  a long optimization time horizon is performed. 
Unfortunately, these pre-computations are strongly dependent on the estimated 
initial position of the targets. Therefore, there is no hope of being 
able to run it in real time taking into account the regularly updated outputs 
from the TMA algorithm. In order to circumvent this difficulty, we 
chose to divide the originally long-term optimization problem into a 
sequence of short-term ones. This procedure is of course sub-optimal, 
but it is compatible with real time  applications and it gives fairly good 
results. 

\subsection*{Horizon splitting technique}
To solve numerically the dynamic programming equations of a finite horizon 
Markov decision process, the parameter horizon $N$ is decisive from the 
computation time and memory consumption points of view. Indeed, 
we need to calculate 
a priori the optimal control associated with each quantization grid point 
for each time step. 
In addition one needs to build and optimize the quantization grids which get larger as the time horizon increases.

In terms of memory, we need to store the transition matrices of size 
$M\times M \times N$, where $M$  is the number of points in the grids 
$\hat{X}^W_t$. In terms of complexity of the algorithm, the computation 
time increases geometrically with $N$, because the size of the 
control space increases geometrically. Unlike scenario 1 and 2 where the 
available controls are in one dimension, we suppose that the control space 
is in three dimensions $(L_1>0, ~L_2>0, ~L_3>0)$, therefore a long horizon 
$N$ will lead to an unsolvable problem. In addition, we have to calculate 
the sound propagation diagram (i.e. cost function $c(l,j)$ in 
algorithm \ref{algo_programmation_dynamique}) for each control point.

To make the computations tractable, we propose a splitting technique: we divide the long-term horizon $N$ 
into $n$ subintervals, and we build the short-term quantization grids once for each subinterval and
apply the short-term dynamic programming 
algorithm on each subinterval, in an iterative manner. Although this is 
sub-optimal compared the original long-term horizon problem, it greatly 
shortens the computation time and makes it possible to run our procedure on line.

 At each time step $t$, 
 two calculations are performed successively as follows.
 First, the TMA process with an Unscented 
 Kalman filter (UKF) is performed to estimate the actual position 
 $\bar{w_t}$ of the target at time $t$. Then one selects the best 
 action given the carrier and target current positions. Suppose 
that $t_0$ is beginning time of sub-optimal interval,  $t_H$ is end time of 
subinterval, and the state of the target at time $t_0$ is estimated by a 
probability distribution $\mathcal{N}(\bar{w}_{t_0}, \Sigma_{t_0})$, which is the 
output of the UKF filter at end time of precedent sub-optimal interval. 
This distribution is used as initial condition to start quantization. Algorithm 
\ref{algo_programmation_dynamique} is applied to pre-calculate optimal control 
for every quantization point. As the position and velocity of the target are unknown, 
(only  their estimations are available), we then search in the 
cloud of quantization grid points the point that is closest to the filtered position and the optimal control corresponding 
to this grid point is applied at this time step. Algorithm \ref{algo2_controle_optimale_en_ligne} is a amended version of algorithm \ref{algo_controle_optimale_en_ligne}. 

\begin{algorithm}[ht]
\SetKwInOut{Input}{input}
\SetKwInOut{Output}{output}
\SetKwInOut{Return}{return}
\caption{Optimal control with known target} 
\label{algo_controle_optimale_en_ligne}
\Input{Optimal control $a_t^\ast$\\
Initial condition of the system $(s_0,w_0)$\\
Grids $\hat{X}_t^W, 0\leq t \leq N$}
\Output{Adaptative control $\hat{a}_t, 0\leq t \leq N$\\
Adaptative $S$ position $s_t, 0\leq t \leq N$}
\Begin{
\For{$t=0$ \KwTo $N-1$}{
Compute estimated target position $\bar{w}_t$ through UKF\\
Find $\hat{w}_t$ the nearest neighbor of $\bar{w}_t$ \\
       among grid $\hat{X}^W_t$\\
$\hat{a}_t = a_t^\ast(s_t,\hat{w}_t)$ optimal control at time $t$
$s_{t+1} = s_t+\hat{a}_t$, update $S$ position
}
}
\end{algorithm}

A consequence of the horizon splitting technique is 
that the quantization grid is automatically updated at the beginning of every 
subinterval. A quantization grid can be considered as an approximation of 
the conditional density given the observations. When the UKF filter has converged, 
the grid points are more concentrated, then accuracy of the MDP 
calculation increases. 

\begin{figure}[ht]
\begin{center} 
\includegraphics[width=0.8\linewidth]{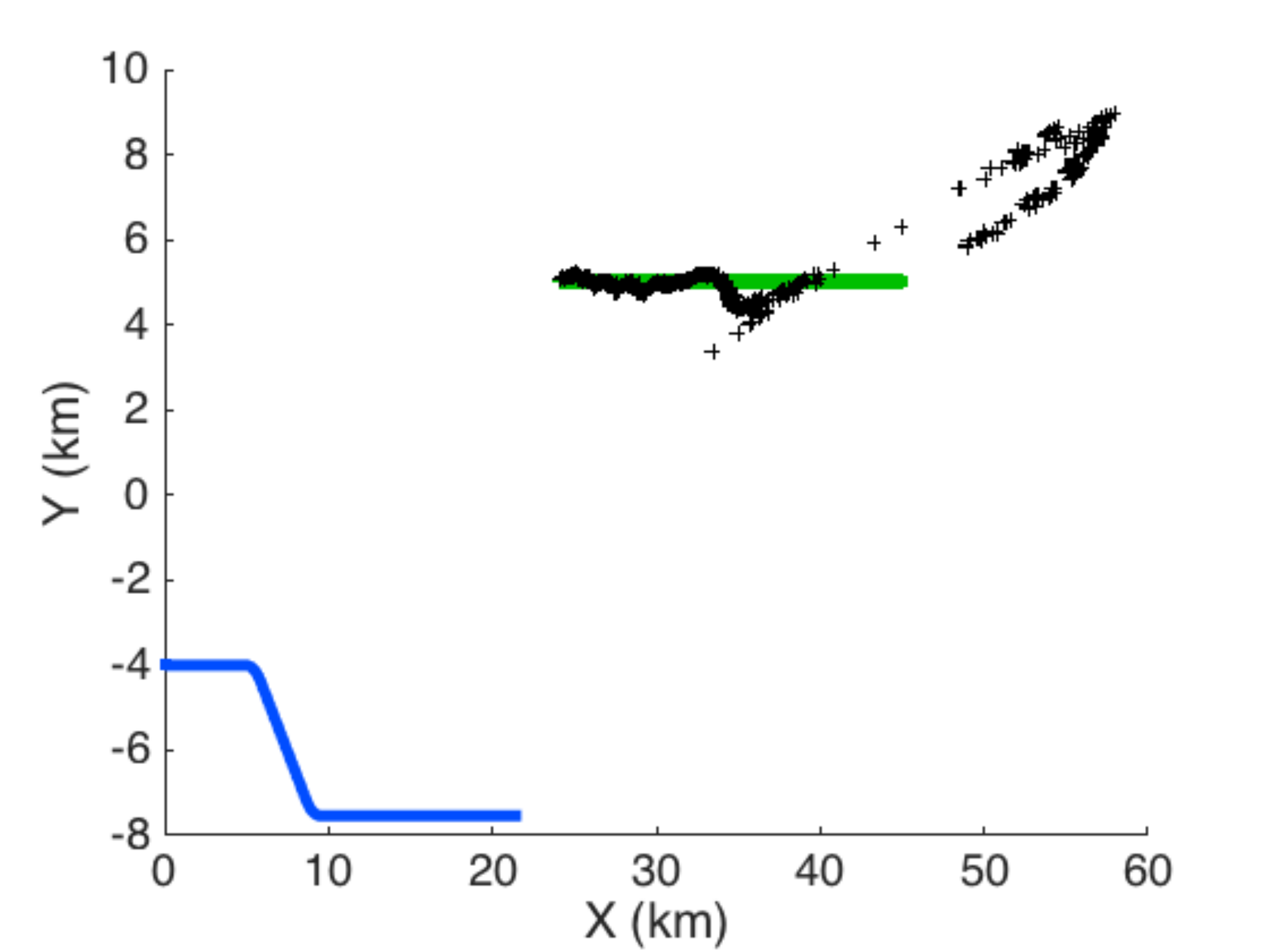}
\caption{Scenario maneuvers and filter without control}
\label{fig_8}
\end{center} 
\end{figure}

\begin{figure}[t]
\begin{center} 
\includegraphics[width=0.8\linewidth]{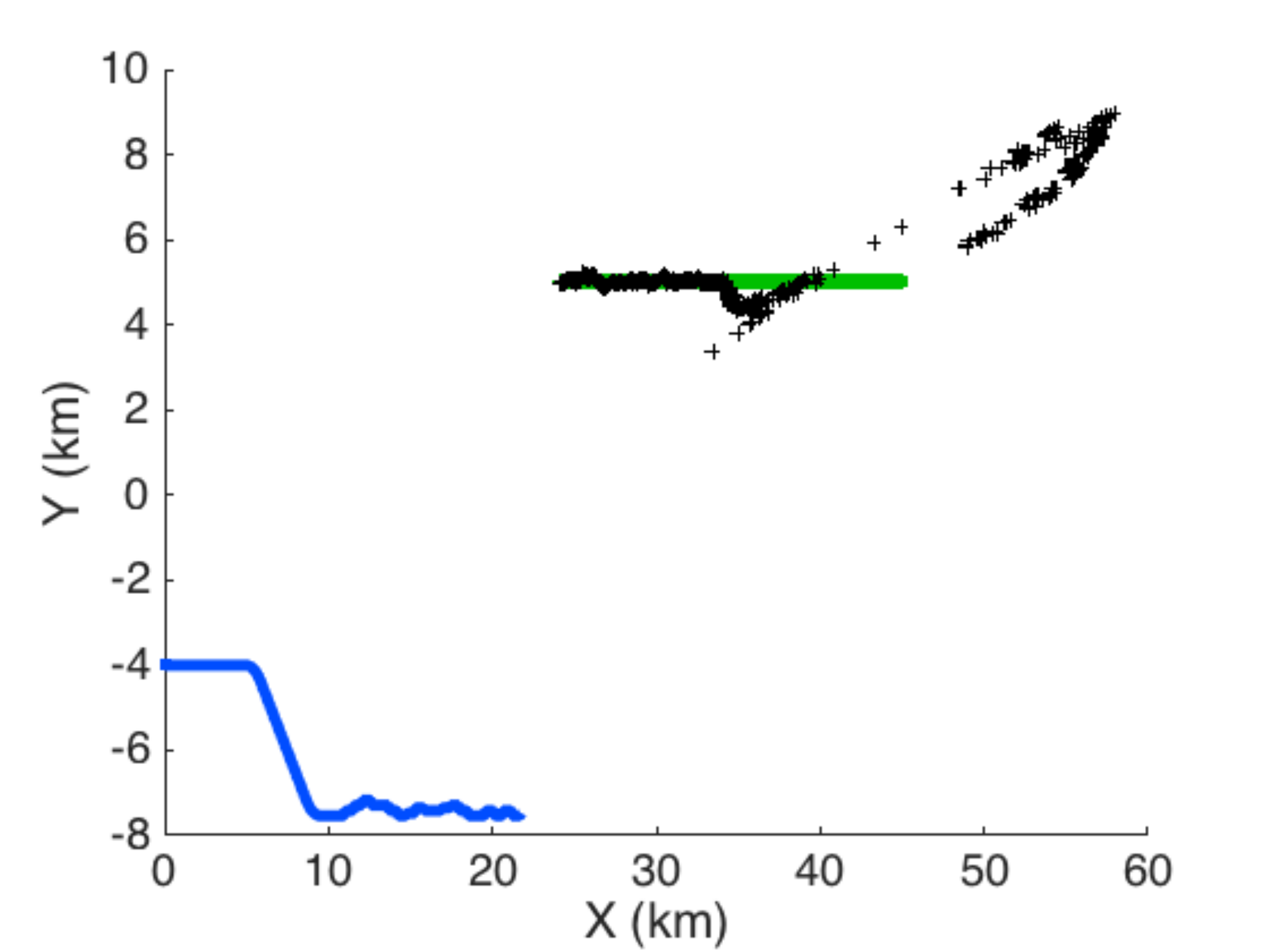}
\caption{Target, submarine and TMA estimation with filtering and control}
\label{fig_9}
\end{center} 
\end{figure}

\begin{figure*}[!t]
\centering
\subfloat{\includegraphics[width=7.0in, height=2.8in] {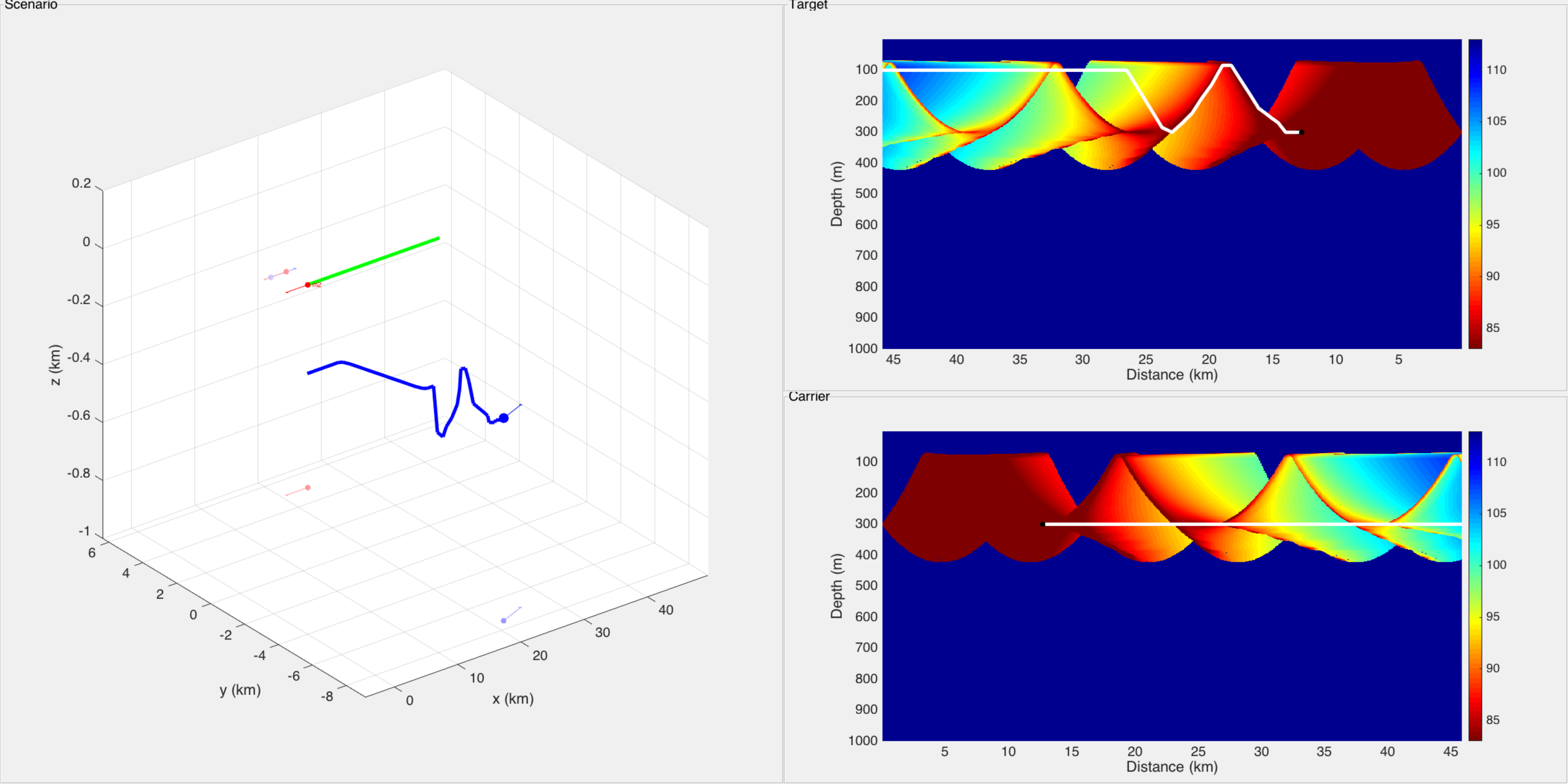}} %
\caption{Target (top right) and carrier (bottom right) sound propagation diagrams and their 3D trajectories (left), with only the target diagram 
under consideration in the cost function. The white curves stand respectively for the carrier depth (top right) and the target depth (bottom right).} 
\label{fig_10}
\end{figure*}

In this scenario there are one carrier submarine and one single target 
in the environment. Objective of the carrier submarine is to detect at best the target. 
The scenario lasts 45 min. The target follows a rectilinear 
motion, and its depth remains constant (300 m). Its position is 
represented by the green curve in Figs.~\ref{fig_8}, 
\ref{fig_9} and 
\ref{fig_12.pdf}, and the green curve in 
Figs.~\ref{fig_10} and 
\ref{fig_13}. 
At each time step, the sensor of the carrier submarine 
processes two measurements (bearing and frequency). Though in this case the 
problem is known to be observable, the performance of the filter is quite poor if 
no maneuvers are scheduled. 
In our scenario, two maneuvers are planned. The first maneuver is a 
right turn at t = 10 min, which lasts 2 minutes, the second is a left turn 
at t = 20 min, which also lasts 2 minutes. These maneuvers allow the filter 
to converge, thereafter the submarine follows a rectilinear motion if no 
control is applied. Fig. \ref{fig_8} illustrates the scenario 
in the horizontal plane without stochastic optimal control (without the depth coordinate as both carrier and target remain at constant depth). 
The trajectory of the carrier submarine $S$ is represented by the blue curve (starting from the left), the real trajectory of the target by the green line (starting from the right) and the black curve 
represents the TMA filter's estimated position of the target.

In order to have an accurate enough estimation of the target position, we have 
divided this scenario into two periods (the time step being $1$ minute).
\begin{itemize} 
\item {\bf Filtering-only period (0-22 min)}. In the first period, only the TMA filter is 
enabled, the period expires at the end of the two maneuvers.
\item {\bf Filtering and optimization period (22-45 min)}. 
This period begins immediately after the first period, 
The three steps algorithms presented in Section \ref{subsect_algorithm} are applied $5$ times for a $5$-minute horizon optimization. 
\end{itemize} 
\vspace{0.1cm}
The numerical results show 
that when the estimation of the target position is of poor quality (error and standard 
deviation too large), efficient control of the trajectory of the submarine is impossible.

As in scenario 1, the cost function is only based on the sound diagram 
emitted by the target. The numerical results are illustrated by two figures. 
Fig.~\ref{fig_9} shows the result of 
filtering. Only the horizontal plane (X-Y) without depth is illustrated. 
Fig.~\ref{fig_10} shows the signal loss due 
to sound propagation emitted by the target (above) and by the carrier 
submarine (below), 
the trajectories in 3D are shown in the left (green: target, blue: carrier). 
The white curve in this diagrams stand respectively for the carrier depth 
and for the target depth. Note that
\begin{itemize} 
\item from $t=22$ min, the trajectory of the carrier submarine is very 
different compared to the original scenario where there is no control (see 
Fig. \ref{fig_8} for comparison). It was modified by the dynamic programming 
procedure, and the quality of the TMA filter is not degraded;
\item the target diagram signal loss shows that the dynamic programming 
has well performed as a controller because, in the second period, the carrier submarine 
remains in the zone of high detection level (reddish zone), despite the fact 
that we replaced the true position of the target by its filtered estimation; 
\item the carrier diagram signal loss is also presented in 
Fig \ref{fig_10}, but it does not 
reflect 
the reality of the past trajectory, because the 
submarine depth was changed after $t = 22$ min leading to a change of its diagram. As the diagram is changing over time, we cannot 
represent it on a static figure. See the web site \cite{zhang16a} for a video 
version. It can be seen that the carrier submarine $S$ can sometimes be easily detected by the 
target because the carrier submarine sound propagation diagram was not taken into account in this scenario. 
\end{itemize} 

\subsection{Scenario 4}
\label{subsect_scenario_4}
We now take fully into account the original mission where two seemingly 
contradictory objectives are to be achieved. The first objective is to pilot the sensor so as to detect at best the 
acoustic signal issued by the target, the second objective is to keep its own detection 
range as low as possible. In this case, the objective function $J$ is a multi criteria aggregation function 
whose optimization will be a trade-off between the conflicting objectives. 
We will take into account the diagrams emitted by the carrier submarine and the target simultaneously in our approach. 

\begin{figure}[t]
\begin{center} 
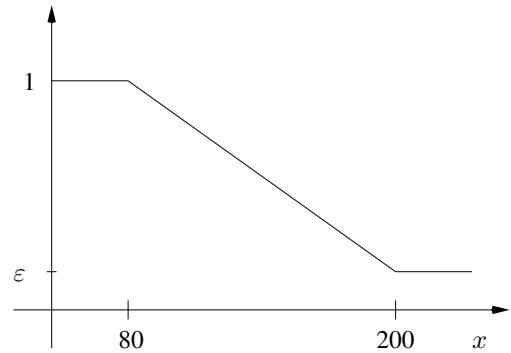
\caption{Multiplicative function $f(x)$}
\label{fig_11}
\end{center} 
\end{figure}

Compared to scenario 3, we use the same scenario parameters, the only 
modification is the cost function $c(l,i)$ used in equation 
(\ref{equ_belleman}). 
Instead of using a weighted sum as in the problem of multi-target
(\ref{equat_fonct_cout_two_targets}), we apply the following method. 
Let $s$ be the position of the carrier submarine, $w$ that of the target, $C^W(s,w)$ 
(resp. $C^S(s,w))$, the value extracted from the sound propagation diagrams emitted by 
the target (resp. by the carrier submarine). These values vary in the interval $[80,200]$. The
objective of the mission is to place the submarine in the area where the 
value of  $C^W(s,w)$ is as low as possible and the value of $C^S(s,w)$
as high as possible. Then the cost function associated to $(s,w)$ is 
defined by
\begin{displaymath} 
  c(s,w) = C^W(s,w)\times f\big(C^S(s,w)\big),
\end{displaymath} 
where $f(x)$ is the following function
\begin{equation} 
\label{fonction_espilon}
f(x) = \left \{
\begin{array} {ll}
1 & ~\mbox{if}~ x<80, \\
ax+b & ~\mbox{if}~ x \in [80,200],\\
\varepsilon & ~\mbox{if}~ x \geq 200,
\end{array} 
\right .
\end{equation} 
where $a$ and $b$ are chosen to satisfy the continuity conditions: $f(80) = 1$ and $f(200) = \varepsilon $. Fig. \ref{fig_11} illustrates the function.
\begin{figure}[h]
\begin{center} 
\includegraphics[width=0.8\linewidth]{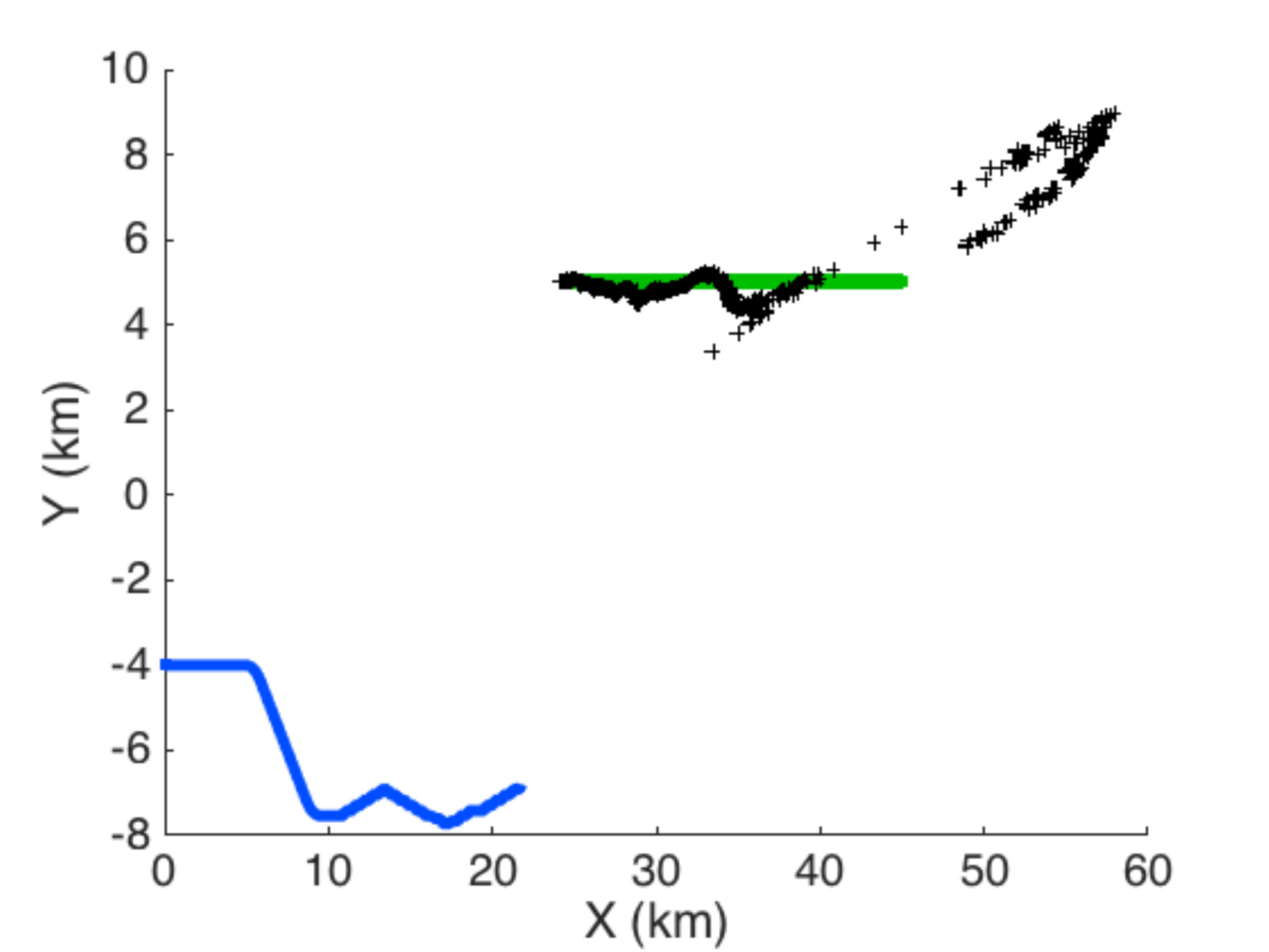}
\caption{Target, submarine and TMA estimation, with consideration of target and submarine diagram in cost function} 
\label{fig_12.pdf}
\end{center} 
\end{figure}
\begin{figure*}[!t]
\centering
\subfloat{\includegraphics[width=7.0in, height=2.8in] {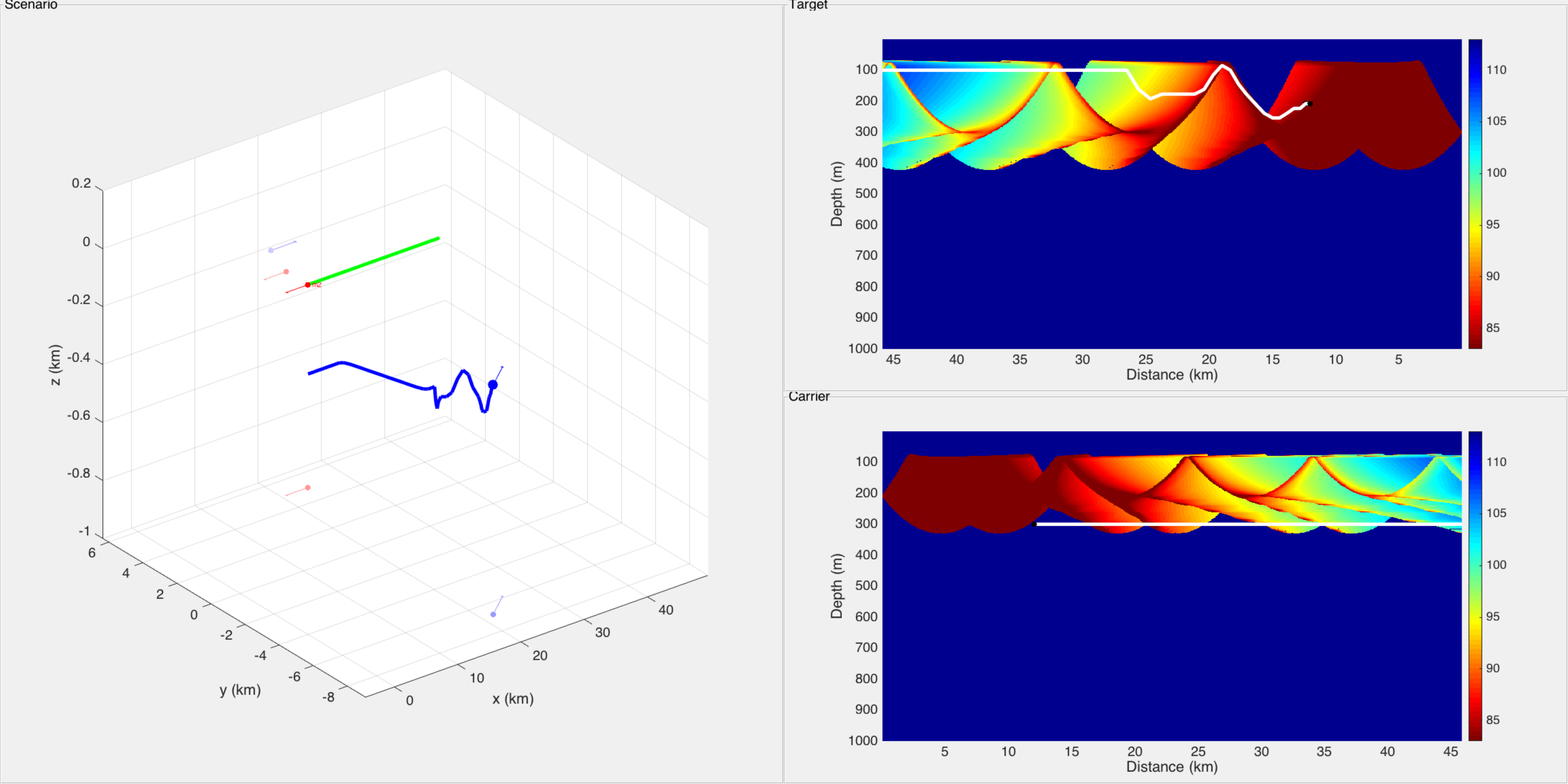}} %
\caption{Target (top right) and carrier (bottom right) sound propagation diagrams and their 3D trajectories (left), with consideration of both target and submarine loss diagrams in the cost function.  The white curves stand respectively for the carrier depth (top right) and the target depth (bottom right).} 
\label{fig_13}
\end{figure*} 

The optimization problem is a minimization of the cost function.
The range of $c(s,w)$ varies in interval $[80\varepsilon, 200]$. When 
it is in the area where $C^S(s,w) \leq 80$ the multiplier is 1, penalizing this area; when 
$C^S(s,w)\geq 200$ the multiplier is $\varepsilon$ making this area of lower cost. In between,
$f$ is a linear function.  Thus we penalize the area where $C^S(s,w)$ is near
 80 without excluding them completely.

The numerical result are presented in Figs. \ref{fig_12.pdf} and \ref{fig_13}. The trajectory of the carrier
submarine is different from that of Figs. 
\ref{fig_9} and \ref{fig_10}, because the cost function and the control policy  are 
different.  Fig.~\ref{fig_13}  shows  the target and carrier signal loss diagrams       respectively.
Note that the detection of the target by the carrier submarine is of lower quality than before, in the sense that 
the carrier submarine is not always in the area of darkest red color. However, the diagram of 
the carrier submarine is now much more  interesting: with a few exceptions, the target 
remains in the blue area, which means that the carrier submarine is undetectable by the target. See the video version of this 
result \cite{zhang16a} for details. This numerical example shows that the two objectives can be realized. 

\section{Conclusion}
In this paper we have proposed an original and effective control strategy for sensor management in a context of 
uncertainty of the targets positions. Four scenarios are considered, from  simple academic case to realistic case.
We showed that it is possible to plan a sensor trajectory even for possibly conflicting objectives. 
With partition technique, the execution time is significantly reduced and the approach can be performed in real time.  
Although this partition method is sub-optimal, it gives very satisfactory results, keeping the carrier submarine position 
in the desired areas. 

\end{document}